\documentclass[11pt]{amsart}
\usepackage{amssymb,amsmath,amsthm}
\usepackage[hmargin=3cm,vmargin=3.5cm]{geometry}
\usepackage{graphicx}
\usepackage{amscd}
\usepackage[all, knot]{xy}
\usepackage{epsfig}
\usepackage{psfrag}
\usepackage{lscape}
\newtheorem{theorem}{Theorem}[section]

\newtheorem{definition}[theorem]{Definition}

\title{Tables of quasi-alternating knots with at most 12 crossings}

\author{SLAVIK JABLAN}

\begin{document}
\maketitle

\begin{abstract}
We are giving tables of quasi-alternating knots with $8\le n \le 12$
crossings. As the obstructions for a knot to be quasialternating we
used homology thickness with regards to Khovanov homology, odd
homology, and Heegaard-Floer homology $\widehat{HFK}$. Except knots
which are homology thick, so cannot be quasialternating, by using
the results of our computations [JaSa1], for one of knots which are
homology thin, knot $11n_{50}$, J.~Greene proved that it is not
quasi-alternating, so this is the first example of homologically
thin knot which is not quasi-alternating [Gr]. In this paper we
provide a few more candidates for homology thin knots for which the
method used by J.~Greene cannot be used to prove that they are not
quasialternating. All computations were performed by A.
Shumakovitch's program \emph{KhoHo}, the program \emph{Knotscape},
the package \emph{Knot Atlas} by Dror Bar-Natan, and our program
\emph{LinKnot}.

\end{abstract}
\newpage


\section{Introduction}

In this paper we present the tables of quasi-alternating (short QA)
knots with $8\le n\le 12$ crossings  based on the thickness of
Khovanov homology $\overline{Kh}$ \cite{Kh1} and odd Khovanov
homology $\overline{Kh}'$ \cite{OzRaSz} and analyzing their minimal
diagrams. Our motivation stems from the fact that such QA-knot
tables are not available from some other sources, so they can be
used by all knot theorists interested for QA knots. Also, they will
maybe stimulate the attempts to prove that some homology thin knots
(so called "candidates" for homology thin knots which are not QA)
are indeed not QA.

C.~Manolescu and P.~Ozsv\'ath show that both Khovanov homology
$\overline{Kh}$ and Heegaard-Floer homology $\widehat{HFK}$ can be
used to detect links which are not QA. Quasi-alternating links are
Khovanov homologically $\sigma $-thin (over $\mathbb{Z}$) and
$\widehat{HFK}$ homologically $\sigma $-thin (over
$\mathbb{Z}$/2$\mathbb{Z}$) \cite{MaOz}. The same property extends
to odd Khovanov homology $\overline{Kh}'$ \cite{OzRaSz} and we will
use this property in the rest of the paper. A knot or link is called
{\it homologically thin} (without qualification), if it is
simultaneously thin with respect to $\overline{Kh}$,
$\widehat{HFK}$, and $\overline{Kh}'$ [Gr, Def. 1.2].

\begin{definition}
The set $Q$ of quasi-alternating links is the smallest set of
links such that

\begin{itemize}
\item the unknot is in $Q$;
\item if the link $L$ has a diagram
$D$ with a crossing $c$ such that

\begin{enumerate}
\item both smoothings of $c$, $L_0$ and $L_{\infty}$, are in $Q$;
\item det($L$)=det$(L_0)$+det$(L_{\infty})$
\end{enumerate}

\end{itemize}

\noindent then $L$ is in $Q$. We say that a crossing $c$
satisfying the properties above is a quasi-alternating crossing of
the diagram $D$ or that $D$ is quasi-alternating at the crossing
$c$ \cite{OzSz,ChKo1}. \label{def}
\end{definition}

The recursive definition makes it difficult to determine if a knot
is quasi-alternating.  It is a challenge to find candidates for
homologically thin knots that are not QA. For a long time, knots
$9_{46}$ = $3,3,-3$ and $10_{140}$ = $4,3,-3$ have been the main
candidates. However, according to A.~Shumakovitch's computations
\cite{Sh2} they are not QA, since are $\overline{Kh}'$-thick
(meaning, they have torsion in their odd Khovanov homology groups),
although they are both $\widehat{HFK}$ and $\overline{Kh}$-thin. By
using the method based on Donaldson's celebrated "Theorem A", which
asserts that the intersection pairing of a smooth, closed, negative
definite 4-manifold is diagonalizable [Don], J.~Greene proved that
homology thin knot $11n_{50}=-2\,2,2\,2,3$ is not QA.

According to Theorem 1 [ChKo1], quasi-alternating links with a
higher number of crossings can be obtained as extension of links
which are already recognized as quasi-alternating \cite{ChKo1,Wi}.

Consider the crossing $c$ in Definition 1 as 2-tangle with marked
endpoints. Using Conway's notation for rational tangles, let
$\varepsilon (c)=\pm 1$, according to whether the overstrand has
positive or negative slope. We will say that a rational 2-tangle
$\tau =C(a_1,\ldots ,a_m)$ extends $c$ if $\tau $ contains $c$ and
$\varepsilon (c)\cdot a_i\ge 1$ for $i=1,\ldots, m$. In particular,
$\tau $ is an alternating rational tangle.

\begin{theorem}
If $L$ is a quasi-alternating link, let $L'$ be obtained by
replacing any quasi-alternating crossing $c$ with an alternating
rational tangle that extends $c$. Then $L'$ is quasi-alternating
$[ChKo1]$.
\end{theorem}

In this paper we give complete computational results for QA knots up
to 12 crossings and the examples of QA knots with at least two
different minimal diagrams, where one is QA and the other is not. We
provide examples of knots and links (short $KL$s) with $n\le 12$
crossings which are homologically thin and have no minimal
quasi-alternating diagrams. In the first version of this paper we
proposed these $KL$s as the candidates for prime homologically thin
links that are not QA, and J.~Greene proved that knot $11n_{50}$ is
not QA [Gr]. Using the method described in his paper it can also be
shown that the link $L11n_{90}$ is not QA, although it is
homologically thin. The remaining candidates given in this paper may
require additional ideas to prove that they are non-QA, if indeed this is the case.
Another attempt to show that some "candidates" are QA is to try to find
their non-minimal diagrams which are QA.  \\

Knots and links are given in Conway notation \cite{Con, JaSa}, which
is implemented in {\it Mathematica} package {\it LinKnot} used for
computations with knots and their distinct diagrams. Odd Khovanov
homology $\overline{Kh}'$ is computed using \emph{KhoHo} by A.
Shumakovitch \cite{Sh1}. All flype-equivalent minimal diagrams of
non-alternating $KL$s up to $n \le 12$ crossings are derived from
alternating link diagrams, using software \emph{LinKnot}. {\it
KnotFind}, the part of the program {\it Knotscape} \cite{HosThi}, is
used for recognition of knots. In addition, we used the criterion
that homologically thick knots are not QA and that for a
quasi-alternating crossing both smoothings must be homologically
thin knots.

\section{Quasi-alternating knots up to 12 crossings}

Since an homologically thick knot cannot be QA, we first selected
knots which are homologically thin. Table 1 gives an overview of the
numbers of non-alternating knots with $8 \le n\le 12$ crossings and
how many among them are homologically thin:
\begin{center}
\begin{table}[h]
\begin{tabular}{|l|c|c|c|c|c|}   \hline
   No. of crossings & 8 & 9 & 10 & 11 & 12 \\ \hline
   No. of non-alternating knots & 3 & 8 & 42 & 185 & 888 \\ \hline
   No. of homologically thin non-alternating knots & 2 & 6 & 31 & 142 & 663 \\ \hline
\end{tabular}
\vspace{0.3cm} \caption{} \label{T1}
\end{table}
\end{center}
Among the knots up to $n=11$ crossings, the following six knots are
$\overline{Kh}$-homologically thin \cite{Sh1} and have a minimal
diagram which is not QA:

\begin{center}
\begin{table}[h]
\begin{tabular}{|c|c|c|c|} \hline
 $9_{46}$  & $3,3,-3$ & $10_{140}$  & $4,3,-3$ \\ \hline
 $11n_{139}$  & $5,3,-3$ & $11n_{107}$ & $-2\,1\,2,3,3$ \\ \hline
 $11n_{50}$  & $-2\,2,2\,2,3$ & $11n_{65}$ & $(3,-2\,1)\,(2\,1,2)$  \\ \hline
\end{tabular}
\vspace{0.3cm} \caption{} \label{T2}
\end{table}
\end{center}

Even columns in Table 2 contain the Conway symbols of these knots.
Four of them, $9_{46}$, $10_{140}$, $11n_{139}$, and $11n_{107}$ are
$\overline{Kh}'$-thick \cite{Sh2}, so they are not QA.

\begin{figure}
\begin{center}
\scalebox{.3}{
\includegraphics{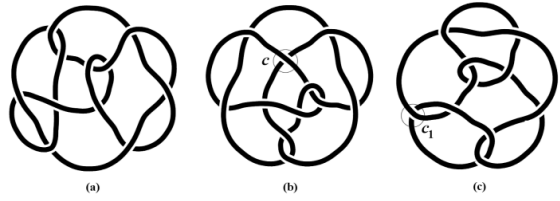}}
\caption{(a) Non-QA minimal diagram $(3,-2\,1)\,(2\,1,2)$ of the
knot $11n_{65}$; (b) the other QA minimal diagram
$6^*2.2\,1.-2\,0.-1.-2$ of the knot $11n_{65}$; (c) QA link
$(2,2+)\,-(2\,1,2)$.}  \label{L11n65}
\end{center} \end{figure}

The knot $11n_{65}$ has two minimal diagrams: $(3,-2\,1)\,(2\,1,2)$
(Fig. \ref{L11n65}a) and $6^*2.2\,1.-2\,0.-1.-2$ (Fig.
\ref{L11n65}b). The first diagram is not QA, and the second is QA.
By smoothing at the crossing $c$, the second diagram resolves into
unknot and QA link $(2,2+)\,-(2\,1,2)$, which resolves into QA knots
$3,2\,1,-2$ and $2\,1\,1,2\,1,-2$ by smoothing the crossing $c_1$
(Fig. \ref{L11n65}c). Moreover, all minimal $KL$ diagrams of the
family derived from knot $11n_{65}$, $(3,-2\,1)\,(p\,1,2)$ and
$6^*2.p\,1.-2\,0.-1.-2$ ($p\ge 2$) which represent the same $KL$
have this property: the first is not QA, and the other is QA.

$12$-crossing knots $12n_{196}$ = $(-3\,1,3)\,(2\,1,2)$ = $6^*2.2
\,0:-3 0.3$, $12n_{393}$ = $8^*2.2\,0:-2\,1\,0$ = $102^*2 0:2\,
0::.-1.-1.-1$ and $12n_{397}$ = $2\,1\,1:-2\,1\,0:2\,0$ have another
minimal diagram which is QA, and the knots given in Table \ref{Lost}
are candidates for homologically thin non-QA knots.

\begin{table}[h]
\begin{tabular}{|c|c|c|c|} \hline
 $12n_{139}$  & $.2.(-2\,1,2).2$ & $12n_{331}$  & $(-3,-2\,-1)\,(3,2+)$ \\ \hline
 $12n_{397}$  & $2\,1\,1:-2\,-1\,0:2\,0$ &   $12n_{414}$  &   $-2\,-1\,0.3.2.2\,0$      \\ \hline
 $12n_{768}$ & $2:-3\,-1\,0:3\,0$ & $12n_{838}$  & $-2.-2.-2\,0.2.2.2\,0$  \\ \hline
\end{tabular}\vspace{0.2cm}
\caption{Candidates for 12-crossing homologically thin non-QA
knots.} \label{Lost}
\end{table}
\begin{figure}
\begin{center}
\scalebox{.4}{
\includegraphics{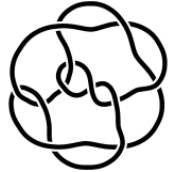}
} \caption{Knot $11n_{50}$. \label{L11n50}}
 \end{center} \end{figure}

\section{Tables of quasi-alternating knots up to 12 crossings}

In the following tables are given QA knots with $8 \le n \le 12$
crossings. Every knot with $8 \le n \le 10$ crossings is given by
its classical symbol and in Conway notation, and every knot with $11
\le n \le 12$ crossings with its DT (Dowker-Thistlethwaite or {\it
Knotscape}) symbol and in Conway notation.

\bigskip

\noindent
\begin{tabular}{|c|c|c|c|} \hline
$8_{20}$ & $3,2\,1,-2$ & $8_{21}$ & $2\,1,2\,1,-2$ \\ \hline
\end{tabular}

\bigskip

\noindent
\begin{tabular}{|c|c|c|c|c|c|} \hline
$9_{43}$ & $2\,1\,1,3,-2$ & $9_{44}$ & $2\, 2,2\,1,-2$ & $9_{45}$ &
$2\,1\,1,2\,1,-2$ \\ \hline

$9_{47}$ & $8^*-2\,0$ & $9_{48}$ & $2\,1,2\,1,-3$ & $9_{49}$ & $-2\,
0:-2\,0:-2\,0$ \\ \hline
\end{tabular}

\small

\bigskip

\noindent
\begin{tabular}{|c|c|c|c|c|c|c|c|} \hline
$10_{125}$ & $5,2\,1,-2$ & $10_{126}$ & $4\,1,3,-2$ & $10_{127}$ & $4\,1,2\,1,-2$ & $10_{129}$ & $3\,2,2\,1,-2$  \\
\hline

$10_{130}$ & $3\,1\,1,3,-2$ & $10_{131}$ & $3\,1\,1,2\,1,-2$ & $10_{133}$ & $2\,3,2\,1,-2$ & $10_{134}$ & $2\,2\,1,3,-2$   \\
\hline

$10_{135}$ & $2\,2\,1,2\,1,-2$ & $10_{137}$ & $2\,2,2\,1\,1,-2$ & $10_{138}$ & $2\,1\,1,2\,1\,1,-2$ &  $10_{141}$ & $4,2\,1,-3$  \\
\hline

$10_{142}$ & $3\,1,3,-2\,-1$ &  $10_{143}$ & $3\,1,3,-3$ &  $10_{144}$ & $3\,1,2\,1,-3$ & $10_{146}$ & $2\,2,2\,1,-3$  \\
\hline

$10_{147}$ & $2\,1\,1,3,-3$ & $10_{148}$ & $(3,2)\,(3,-2)$ & $10_{149}$ & $(3,2)\,(2\,1,-2)$ & $10_{150}$ & $(2\,1,2)\,(3,-2)$   \\
\hline

$10_{151}$ & $(2\,1,2)\,(2\,1,-2)$ & $10_{155}$ & $-3:2:2$ & $10_{156}$ & $-3:2:2 0$ & $10_{157}$ & $-3:2\,0:2\,0$   \\
\hline

 $10_{158}$ & $-3\,0:2:2$ & $10_{159}$ & $-3\,0:2:2\,0$ & $10_{160}$ & $-3\,0:2\,0:2\,0$ & $10_{162}$ & $-3\,0:-2\,0:-2\,0$   \\
\hline

 $10_{163}$ & $8^*-3\,0$ &$10_{164}$ & $8^*2:-2\,0$  & $10_{165}$ & $8^*2:.-2\,0$ & &  \\
\hline
\end{tabular}

\normalsize

\bigskip

\begin{landscape}
\small

\noindent
\begin{tabular}{|c|c|c|c|c|c|c|c|} \hline
$K11n1$ & $2\,3,2\,1\,1,-2$ & $K11n2$ & $2\,2\,1,2\,1\,1,-2$ & $K11n3$ & $2\,2\,1,2\,2,-2$ & $K11n4$ & $(2\,1\,1,2)\,(3,-2)$   \\
\hline

$K11n5$ & $(2\,1\,1,2) (2\,1,-2)$ & $K11n7$ & $(2\,1\,1,-2)\,(2\,1,2)$ & $K11n8$ & $(2\,2,-2) (2\,1,2)$ & $K11n10$ & $(2\,1\,1,-2)\,(3,2)$   \\
\hline

$K11n11$ & $(2\,2,-2)\,(3,2)$ & $K11n13$ & $5,2\,1\,1,-2$ & $K11n14$ & $4\,1,2\,1\,1,-2$ & $K11n15$ & $4\,1,2\,2,-2$   \\
\hline

$K11n16$ & $3\,2,2\,1\,1,-2$ & $K11n17$ & $3\,1\,1,2\,1\,1,-2$ & $K11n18$ & $3\,1\,1,2\,2,-2$ & $K11n21$ & $(2\,1,2+) (2\,1,-2)$   \\
\hline

$K11n22$ & $(2\,1,2+1)\,-(2\,1,2)$ & $K11n23$ & $(2\,1,2+)\,-(3,2)$ & $K11n25$ & $(3,2+)\,(2\,1,-2)$ & $K11n26$ & $(3,2+)\,-(2\,1,2)$   \\
\hline

$K11n28$ & $2\,4,2\,1,-2$ & $K11n29$ & $2\,3\,1,2\,1,-2$ & $K11n30$ & $2\,3\,1,3,-2$ & $K11n32$ & $(2\,2,2)\,(2\,1,-2)$   \\
\hline

$K11n33$ & $(2\,2,2)\,(3,-2)$ & $K11n35$ & $.(2\,1,-2).2\,0$ & $K11n36$ & $.(3,-2).2\,0$ & $K11n37$ & $(3,2+)\,(3,-2)$   \\
\hline

$K11n40$ & $.(2\,1,-2).2$ & $K11n41$ & $.(3,-2).2$ & $K11n43$ & $.2\,0.(2\,1,-2)$ & $K11n44$ & $.2\,0.(3,-2)$   \\
\hline

$K11n48$ & $3\,1,2\,2,-2\,-1$ & $K11n51$ & $2\,1\,3,2\,1,-2$ & $K11n52$ & $2\,1\,2\,1,2\,1,-2$ & $K11n53$ & $2\,1\,2\,1,3,-2$   \\
\hline

$K11n54$ & $2\,1\,1\,2,2\,1,-2$ & $K11n55$ & $2\,1\,1\,1\,1,2\,1,-2$ & $K11n56$ & $2\,1\,1\,1\,1,3,-2$ & $K11n58$ & $3\,1\,2,2\,1,-2$   \\
\hline

$K11n59$ & $3\,1\,1\,1,2\,1,-2$ & $K11n60$ & $3\,1\,1\,1,3,-2$ & $K11n62$ & $4\,2,2\,1,-2$ & $K11n63$ & $4\,1\,1,2\,1,-2$   \\
\hline

$K11n64$ & $4\,1\,1,3,-2$ & $K11n65$ & $(3,-2\,-1)\,(2\,1,2)$ & $K11n66$ & $(2\,1,-3)\,(2\,1,2)$ & $K11n68$ & $(3,2\,1)\,(2\,1,-2)$   \\
\hline

$K11n69$ & $(3,2\,1)\,(3,-2)$ & $K11n71$ & $2\,1,3,2\,1,-2$ & $K11n72$ & $2\,1,2\,1,2\,1,-2$ & $K11n75$ & $3,2\,1,2\,1,-2$   \\
\hline

$K11n76$ & $3,3,2\,1,-2$ & $K11n78$ & $3,2\,1,3,-2$ & $K11n82$ & $4,2\,2,-3$ & $K11n83$ & $3\,1,2\,2,-3$   \\
\hline

$K11n84$ & $2\,2,2\,2,-3$ & $K11n85$ & $2\,1\,1\,1,3,-3$ & $K11n86$ & $-2\,-1\,-1:2:2\,0$ & $K11n87$ & $2\,1\,2,2\,1,-3$   \\
\hline

$K11n89$ & $3\,1,2\,1\,1,-3$ & $K11n90$ & $3\,1,2\,1\,1,-2\,-1$ & $K11n91$ & $4,2\,1\,1,-3$ & $K11n93$ & $-2\,-1\,0:2\,1:2$   \\
\hline

$K11n94$ & $8^*-2\,-1\,-1\,0$ & $K11n95$ & $2\,0.-2\,-1.-2\,0.2$ & $K11n98$ & $(2\,1,-3)\,(3,2)$ & $K11n99$ & $(3,-2\,-1)\,(3,2)$   \\
\hline

$K11n100$ & $2\,2\,1,3,-3$ & $K11n101$ & $2\,3,2\,1,-3$ & $K11n103$ & $2\,1\,1,2\,1\,1,-3$ & $K11n105$ & $2\,1\,1\,1,2\,1,-3$  \\
\hline

$K11n106$ & $2\,1\,2,3,-3$ & $K11n108$ & $-3\,0:2\,1:2$ & $K11n109$ & $-2\,-1\,0:-3\,0:-2\,0$ & $K11n110$ & $-2\,-1\,-1\,0:2:2\,0$  \\
\hline

$K11n112$ & $-2\,-1\,-1\,0:2:2$ & $K11n113$ & $-2\,-2:2:2$ & $K11n114$ & $-4:2:2\,0$ & $K11n115$ & $-2\,-2:2\,0:2\,0$  \\
\hline

$K11n117$ & $-4\,0:2\,0:2\,0$ & $K11n118$ & $2\,0.-3.-2\,0.2$ & $K11n119$ & $8^*2\,1::-2\,0$ & $K11n120$ & $8^*-2\,-1\,0:.2\,0$   \\
\hline

$K11n121$ & $-2\,-2:-2\,0:-2\,0$ & $K11n122$ & $3\,2,3,-3$ & $K11n123$ & $2.2\,1.-2.2$ & $K11n124$ & $-2\,-2\,0:2:2\,0$   \\
\hline

$K11n125$ & $-2\,-1\,-1:2\,0:2\,0$ & $K11n127$ & $-2\,-1\,-1:-2\,0:-2\,0$ & $K11n128$ & $8^*-2\,-1\,0:.2$ & $K11n129$ & $-2\,-1\,0:3:2$   \\
\hline

$K11n130$ & $8^*-2\,-1\,0:2$ & $K11n131$ & $-3\,0:2\,1\,0:2$ & $K11n132$ & $-2\,-1\,0:3\,0:2$ & $K11n134$ & $2.-2\,-1.2.2$  \\
\hline

$K11n136$ & $2\,-1\,0:-2\,-1\,0:-2\,0$ & $K11n137$ & $3\,1\,1,2\,1,-3$ & $K11n140$ & $4\,1,2\,1,-3$ & $K11n141$ & $4\,1,3,-2\,-1$  \\
\hline

$K11n142$ & $2.3.-2.2\,0$ & $K11n144$ & $-2\,-1\,-1\,0:-2\,0:-2\,0$ & $K11n146$ & $8^*2\,1\,0:-2\,0$ & $K11n147$ & $8^*-2\,-1\,0:2\,0$  \\
\hline

$K11n148$ & $9^*.-2\,-1$ & $K11n149$ & $8^*-3\,-1\,0$ & $K11n150$ & $8^*-4\,0$ & $K11n153$ & $8^*3::-2\,0$   \\
\hline

$K11n154$ & $8^*-3\,0:.2\,0$ & $K11n155$ & $8^*3\,0::-2\,0$ & $K11n156$ & $10^*-2\,0$ & $K11n157$ & $8^*2.-2\,0.2\,0$  \\
\hline

$K11n158$ & $8^*3\,0:-2\,0$ & $K11n159$ & $8^*-3\,0:2\,0$ & $K11n160$ & $8^*2:2:-2\,0$ & $K11n161$ & $8^*3:-2\,0$  \\
\hline

$K11n162$ & $8^*2:.-2\,0:.2$ & $K11n163$ & $8^*-3\,0.2\,0$ & $K11n164$ & $2.2.-2.2\,0.2\,0$ & $K11n165$ & $8^*2.2\,0.-2\,0$  \\
\hline

$K11n166$ & $8^*2:2\,0:-2\,0$ & $K11n167$ & $2.2.-2.2.2\,0$ & $K11n168$ & $8^*-3\,0::2\,0$ & $K11n169$ & $-3\,-1\,0:2:2$  \\
\hline

$K11n170$ & $2\,0.3.-2.2$ & $K11n171$ & $-4\,0:2:2$ & $K11n172$ & $8^*2\,0:2\,0:-2\,0$ & $K11n173$ & $8^*2\,0:-2\,0:2\,0$  \\
\hline

$K11n174$ & $9^*2.-2$ & $K11n175$ & $8^*2:.2:.-2\,0$ & $K11n176$ & $8^*2:-2\,0:2\,0$ & $K11n177$ & $9^*2\,0.-2$  \\
\hline

$K11n178$ & $10^{**}-2\,0$ & $K11n179$ & $2.-3.2.2\,0$ & $K11n180$ & $-3\,-1\,0:-2\,0:-2\,0$ & $K11n181$ & $-4\,0:-2\,0:-2\,0$  \\
\hline

$K11n182$ & $9^*.2\,0:.-2$ & $K11n184$ & $9^*.2:.-2$ & $K11n185$ & $9^*.-3$ & & \\
\hline

\end{tabular}

\end{landscape} \normalsize

\tiny \noindent
\begin{tabular}{|c|c|c|c|c|c|} \hline
$K12n1$ & $(-2\,-1\,-1,2)\,1\,(2\,1,2)$ & $K12n2$ &
$(2\,1\,1,-2)\,1\,(2\,1,2)$ & $K12n3$ &
$(-2\,-1\,-1,-2)\,1\,(2\,1,2)$ \\ \hline $K12n4$ &
$(2\,1\,1,2)\,1\,(-2\,-1,-2)$ & $K12n5$ &
$(2\,1\,1,2)\,1\,(2\,1,-2)$ & $K12n6$ & $(2\,1\,1,2)\,1\,(-2\,-1,2)$
\\ \hline $K12n8$ & $(2\,1\,1,-2)\,1\,(3,2)$ & $K12n9$ &
$(-2\,-1\,-1,-2)\,1\,(3,2)$ & $K12n11$ & $2\,3\,1,-2\,-1\,-1,2$ \\
\hline $K12n12$ & $2\,3\,1,2\,1\,1,-2$ & $K12n13$ &
$-2\,-3\,-1,2\,1\,1,2$ & $K12n14$ & $(-2\,-1\,-1,2)\,(2\,1\,1,2)$ \\
\hline $K12n15$ & $(2\,1\,1,-2)\,(2\,1\,1,2)$ & $K12n17$ &
$(2\,2,2)\,(-2\,-1\,-1,2)$ & $K12n18$ & $(2\,2,2)\,(2\,1\,1,-2)$ \\
\hline $K12n21$ & $.2.(-2\,-1\,-1,2)\,0$ & $K12n22$ &
$.2.(2\,1\,1,-2)\,0$ & $K12n24$ & $(-2\,-1\,-1,2)\,(3,2+1)$ \\
\hline $K12n26$ & $.2.(2,-2\,-1\,-1)$ & $K12n27$ & $.2.(-2,2\,1\,1)$
& $K12n29$ & $.2.(2,-2\,-1\,-1)\,0$ \\ \hline $K12n30$ &
$.2.(-2,2\,1\,1)\,0$ & $K12n32$ & $.2.(-2\,-1\,-1,2)$ & $K12n33$ &
$.2.(2\,1\,1,-2)$ \\ \hline $K12n35$ & $2\,1\,2\,1,-2\,-1\,-1,2$ &
$K12n36$ & $2\,1\,2\,1,2\,1\,1,-2$ & $K12n37$ &
$-2\,-1\,-2\,-1,2\,1\,1,2$ \\ \hline $K12n38$ &
$2\,1\,1\,1\,1,-2\,-1\,-1,2$ & $K12n39$ & $2\,1\,1\,1\,1,2\,1\,1,-2$
& $K12n40$ & $-2\,-1\,-1\,-1\,-1,2\,1\,1,2$ \\ \hline $K12n42$ &
$3\,1\,1\,1,-2\,-1\,-1,2$ & $K12n43$ & $3\,1\,1\,1,2\,1\,1,-2$ &
$K12n44$ & $-3\,-1\,-1\,-1,2\,1\,1,2$ \\ \hline $K12n46$ &
$4\,1\,1,-2\,-1\,-1,2$ & $K12n47$ & $4\,1\,1,2\,1\,1,-2$ & $K12n48$
& $-4\,-1\,-1,2\,1\,1,2$ \\ \hline $K12n49$ &
$(-2\,-1\,-1,2)\,(3,2\,1)$ & $K12n50$ & $(2\,1\,1,-2)\,(3,2\,1)$ &
$K12n52$ & $(2\,1\,1,2)\,(-3,2\,1)$ \\ \hline $K12n53$ &
$(2\,1\,1,2)\,(3,-2\,-1)$ & $K12n55$ & $-2\,-1\,-1,2\,1,2\,1,2$ &
$K12n58$ & $2\,1\,1,2\,1,2\,1,-2$ \\ \hline $K12n60$ &
$2\,1\,1,-2\,-1,2\,1,2$ & $K12n61$ & $2\,1\,1,2\,1,-2\,-1,2$ &
$K12n62$ & $-2\,-1\,-1,2,3,2\,1$ \\ \hline $K12n64$ &
$2\,1\,1,2,3,-2\,-1$ & $K12n66$ & $-2\,-1\,-1,3,2\,1,2$ & $K12n69$ &
$(2\,2,2)\,1\,(2\,1,-2)$ \\ \hline $K12n70$ &
$(2\,2,2)\,1\,(-2\,-1,-2)$ & $K12n71$ & $(2\,1,2+1)\,1\,(-2\,-1,2)$
& $K12n72$ & $(2\,1,2+1)\,1\,(2\,1,-2)$ \\ \hline $K12n73$ &
$(2\,1,2+1)\,1\,(-2\,-1,-2)$ & $K12n74$ & $(3,2+1)\,1\,(-2\,-1,2)$ &
$K12n75$ & $(3,2+1)\,1\,(2\,1,-2)$ \\ \hline $K12n76$ &
$(3,2+1)\,1\,(-2\,-1,-2)$ & $K12n77$ & $2\,4\,1,-2\,-1,2$ & $K12n78$
& $2\,4\,1,2\,1,-2$ \\ \hline $K12n79$ & $-2\,-4\,-1,2\,1,2$ &
$K12n80$ & $(2\,2\,1,2)\,(-2\,-1,2)$ & $K12n81$ &
$(2\,2\,1,2)\,(2\,1,-2)$ \\ \hline $K12n83$ &
$(2\,2\,1,-2)\,(2\,1,2)$ & $K12n84$ & $(-2\,-2\,-1,2)\,(2\,1,2)$ &
$K12n85$ & $(-2\,-1,2),2,(2\,1,2)$ \\ \hline $K12n86$ &
$(2\,1,-2),2,(2\,1,2)$ & $K12n88$ & $(2\,1,2),-2,(2\,1,2)$ &
$K12n90$ & $(3,2),2,(2\,1,-2)$ \\ \hline $K12n92$ &
$(3,2),-2,(2\,1,2)$ & $K12n94$ & $(2\,3,2)\,(-2\,-1,2)$ & $K12n95$ &
$(2\,3,2)\,(2\,1,-2)$ \\ \hline $K12n97$ & $(2\,2,2+1)\,(-2\,-1,2)$
& $K12n98$ & $6^*(-2\,-1,2).2\,1$ & $K12n99$ & $6^*(2\,1,-2).2\,1$ \\
\hline $K12n100$ & $.2.(-2\,-1,2)\,0.2$ & $K12n101$ &
$.2.(2\,1,-2)\,0.2$ & $K12n103$ & $(4\,1,2)\,(-2\,-1,2)$ \\ \hline
$K12n104$ & $(4\,1,2)\,(2\,1,-2)$ & $K12n106$ &
$(4\,1,-2)\,(2\,1,2)$ & $K12n107$ & $(-4\,-1,2)\,(2\,1,2)$ \\ \hline
$K12n108$ & $(3\,1\,1,2)\,(-2\,-1,2)$ & $K12n109$ &
$(3\,1\,1,2)\,(2\,1,-2)$ & $K12n111$ & $(3\,1\,1,-2)\,(2\,1,2)$ \\
\hline $K12n112$ & $(-3\,-1\,-1,2)\,(2\,1,2)$ & $K12n113$ &
$(5,2)\,(-2\,-1,2)$ & $K12n114$ & $(5,2)\,(2\,1,-2)$ \\ \hline
$K12n116$ & $(3\,2,2)\,(-2\,-1,2)$ & $K12n117$ &
$(3\,2,2)\,(2\,1,-2)$ & $K12n122$ & $6^*(-2\,-1,2).2\,1\,0$ \\
\hline $K12n123$ & $6^*(2\,1,-2).2\,1\,0$ & $K12n125$ &
$6^*(2,-2\,-1).2\,1$ & $K12n126$ & $6^*(-2,2\,1).2\,1$ \\ \hline
$K12n127$ & $6^*(2,-2\,-1).2\,1\,0$ & $K12n128$ &
$6^*(-2,2\,1).2\,1\,0$ & $K12n130$ & $(-2\,-1,2),2,(2,2\,1)$ \\
\hline $K12n131$ & $(2\,1,-2),2,(2,2\,1)$ & $K12n133$ &
$(2\,1,2),-2,(2,2\,1)$ & $K12n135$ & $(3,2),2,(-2,2\,1)$ \\ \hline
$K12n137$ & $(3,2),-2,(2,2\,1)$ & $K12n140$ & $.2.(2\,1,-2).2$ &
$K12n142$ & $.2\,1\,1.-2\,-1.2\,0$ \\ \hline $K12n144$ &
$.2\,1.-2\,-1.3$ &
$K12n145$ & $.2\,1.-2\,-1.3\,0$ & $K12n146$ & $.2\,2.-2\,-1.2\,0$ \\
\hline $K12n147$ & $.2\,1.-2\,-1.2\,1$ & $K12n150$ &
$2\,1\,1\,1\,1\,1,-2\,-1,2$ & $K12n151$ & $2\,1\,1\,1\,1\,1,2\,1,-2$
\\ \hline $K12n152$ & $2\,1\,1\,1\,1\,1,-2\,-1,-2$ & $K12n154$ &
$2\,2\,3,2\,1,-2$ & $K12n155$ & $2\,2\,3,-2\,-1,-2$ \\ \hline
$K12n157$ & $(2\,1,2\,1+1)\,(2\,1,-2)$ & $K12n158$ &
$(2\,1,2\,1+1)\,(-2\,-1,-2)$ & $K12n159$ & $2\,2\,2\,1,-2\,-1,2$ \\
\hline $K12n160$ & $2\,2\,1\,1\,1,-2\,-1,2$ & $K12n161$ &
$2\,2\,1\,1\,1,2\,1,-2$ & $K12n162$ & $2\,2\,1\,1\,1,-2\,-1,-2$ \\
\hline $K12n163$ & $2\,1\,2\,1\,1,-2\,-1,2$ & $K12n164$ &
$2\,1\,2\,1\,1,2\,1,-2$ & $K12n165$ & $2\,1\,2\,1\,1,-2\,-1,-2$ \\
\hline $K12n166$ & $4\,2\,1,-2\,-1,2$ & $K12n167$ &
$4\,2\,1,2\,1,-2$ & $K12n168$ & $4\,2\,1,-2\,-1,-2$ \\ \hline
$K12n169$ & $3\,1\,2\,1,-2\,-1,2$ & $K12n170$ & $3\,1\,2\,1,2\,1,-2$
& $K12n171$ & $3\,1\,2\,1,-2\,-1,-2$ \\ \hline $K12n173$ &
$8^*(-2\,-1,2)$ & $K12n174$ & $8^*(2\,1,-2)$ & $K12n176$ &
$.(2\,1,2).-2.2\,0$ \\ \hline $K12n177$ & $(2\,2,3)\,(-2\,-1,2)$ &
$K12n178$ & $(2\,2,3)\,(2\,1,-2)$ & $K12n180$ &
$(2\,2,-3)\,(2\,1,2)$ \\ \hline $K12n181$ & $(-2\,-2,3)\,(2\,1,2)$ &
$K12n182$ & $(2\,1\,1,2\,1)\,(-2\,-1,2)$ & $K12n183$ &
$(2\,1\,1,2\,1)\,(2\,1,-2)$ \\ \hline $K12n185$ &
$(3\,1,2\,1)\,(-2\,-1,2)$ & $K12n186$ & $(3\,1,2\,1)\,(2\,1,-2)$ &
$K12n188$ & $(-3\,-1,2\,1)\,(2\,1,2)$ \\ \hline $K12n189$ &
$(3\,1,-2\,-1)\,(2\,1,2)$ & $K12n190$ & $(4,3)\,(-2\,-1,2)$ &
$K12n191$ & $(4,3)\,(2\,1,-2)$ \\ \hline $K12n193$ &
$(3\,1,3)\,(-2\,-1,2)$ & $K12n194$ & $(3\,1,3)\,(2\,1,-2)$ &
$K12n196$ & $(-3\,-1,3)\,(2\,1,2)$ \\ \hline $K12n197$ &
$(3\,1,-3)\,(2\,1,2)$ & $K12n201$ & $(-2\,-1,2):2:2$ & $K12n202$ &
$(2\,1,-2):2:2$ \\ \hline $K12n204$ & $.2.-2.(2,2\,1)\,0$ &
$K12n205$ & $.(-2\,-1,2).2.2\,0$ & $K12n206$ & $.(2\,1,-2).2.2\,0$
\\ \hline $K12n208$ & $.2.2.(2,-2\,-1)\,0$ & $K12n209$ &
$.2.2.(-2,2\,1)\,0$ & $K12n211$ & $(2\,1,2\,1)\,1\,(-2\,-1,2)$ \\
\hline $K12n212$ & $.2.2.(-2\,-1,2)\,0$ & $K12n213$ &
$.2.2.(2\,1,-2)\,0$ & $K12n215$ & $(4,2\,1)\,(-2\,-1,2)$ \\ \hline
$K12n216$ & $(4,2\,1)\,(2\,1,-2)$ & $K12n219$ &
$2\,1\,1,-2\,-1,2,2\,1$ & $K12n222$ & $2\,1\,1,2\,1,-2,2\,1$ \\
\hline $K12n223$ & $-2\,-1\,-1,2\,1,2,2\,1$ & $K12n224$ &
$2\,2,-2\,-1,2,2\,1$ & $K12n226$ & $.(2,-2\,-1).2.2\,0$ \\ \hline
$K12n227$ & $.(-2,2\,1).2.2\,0$ & $K12n233$ & $6\,1,-2\,-1,2$ &
$K12n234$ & $6\,1,2\,1,-2$ \\ \hline $K12n235$ & $6\,1,-2\,-1,-2$ &
$K12n236$ & $3\,3\,1,-2\,-1,2$ & $K12n237$ & $3\,3\,1,2\,1,-2$ \\
\hline $K12n238$ & $3\,3\,1,-2\,-1,-2$ & $K12n239$ &
$3\,2\,1\,1,-2\,-1,2$ & $K12n240$ & $3\,2\,1\,1,2\,1,-2$ \\ \hline
$K12n241$ & $3\,2\,1\,1,-2\,-1,-2$ & $K12n245$ &
$(3,3+1)\,(-2\,-1,2)$ & $K12n246$ & $(3,3+1)\,(2\,1,-2)$ \\ \hline
$K12n247$ & $(3,3+1)\,(-2\,-1,-2)$ & $K12n248$ & $5\,1\,1,-2\,-1,2$
& $K12n249$ & $5\,1\,1,2\,1,-2$ \\ \hline $K12n250$ &
$5\,1\,1,-2\,-1,-2$ & $K12n252$ & $2:(2,-2\,-1)\,0:2\,0$ & $K12n253$
& $2:(2,-2\,-1)\,0:-2\,0$ \\ \hline $K12n254$ &
$2:(-2,2\,1)\,0:-2\,0$ & $K12n255$ & $2:(-2,2\,1)\,0:2\,0$ &
$K12n259$ & $2:(2,2\,1)\,0:-2\,0$ \\ \hline $K12n261$ &
$2\,1\,1,3,2,-2\,-1$ & $K12n262$ & $2:(-2\,-1,2)\,0:2\,0$ &
$K12n263$ & $2:(2\,1,-2)\,0:2\,0$ \\ \hline $K12n266$ &
$(2\,1,-2):2\,0:2\,0$ & $K12n269$ & $(2\,2,2)\,(-3,2\,1)$ &
$K12n270$ & $(2\,2,2)\,(3,-2\,-1)$ \\ \hline $K12n271$ &
$.2.-2\,-1.2\,1\,1\,0$ & $K12n272$ & $.2\,1\,1.-2\,-1.2$ & $K12n274$
& $.3.-2\,-1.2\,1\,0$ \\ \hline $K12n275$ & $.2.-2\,-1.2\,2\,0$ &
$K12n277$ & $.3\,1.-2\,-1.2$ & $K12n278$ & $.4.-2\,-1.2$ \\ \hline
$K12n280$ & $2.-2\,0.-2.2\,1\,1\,0$ & $K12n281$ &
$2\,1\,1\,0.2.-2.2\,0$ & $K12n283$ & $8^*2\,1\,1.-2\,0$ \\ \hline
$K12n284$ & $8^*-2\,-1\,-1::2\,0$ & $K12n285$ &
$8^*-2\,-1\,-1::-2\,0$ & $K12n286$ & $.2\,1\,1\,1.-2\,0.2$ \\ \hline
$K12n287$ & $8^*2\,1\,1::-2\,0$ & $K12n288$ & $-2\,-2\,-1,2\,2\,1,2$
& $K12n289$ & $2\,2\,1,2\,2\,1,-2$ \\ \hline $K12n290$ &
$(-2\,-2\,-1,2)\,(3,2)$ & $K12n291$ & $(2\,2\,1,-2)\,(3,2)$ &
$K12n294$ & $-3\,-1,2\,1\,1,2\,1\,1$ \\ \hline $K12n295$ &
$-2\,-1\,-2,2\,1\,1,2\,1$ & $K12n296$ & $2\,1\,2,2\,1\,1,-2\,-1$ &
$K12n297$ & $3\,1,-2\,-2,2\,1\,1$ \\ \hline $K12n298$ &
$8^*-2\,-1\,-1\,0:.2\,0$ & $K12n299$ & $-2.-2\,0.-2.2\,1\,1\,0$ &
$K12n300$ & $-2\,-1\,-1\,-1,2\,1\,1,2\,1$ \\ \hline $K12n301$ &
$2\,1\,1\,1,2\,1\,1,-2\,-1$ & $K12n302$ & $.2\,1\,1\,1.-2.2$ &
$K12n303$ & $4\,1,-2\,-2\,-1,2$ \\ \hline $K12n304$ &
$4\,1,2\,2\,1,-2$ & $K12n305$ & $-4\,-1,2\,2\,1,2$ & $K12n306$ &
$3\,1\,1,-2\,-2\,-1,2$ \\ \hline $K12n307$ & $3\,1\,1,2\,2\,1,-2$ &
$K12n308$ & $-3\,-1\,-1,2\,2\,1,2$ & $K12n311$ &
$-2\,-1\,-1\,0.2.2.2\,0$ \\ \hline $K12n312$ &
$-2\,-1\,-1:2:2\,1\,0$ & $K12n315$ & $8^*3.-2\,-1$ & $K12n316$ &
$2\,0.2\,1.-2\,0.3\,0$ \\ \hline $K12n317$ &
$(2\,1\,1,-2\,-1)\,(3,2)$ & $K12n319$ & $2\,0.-2\,-1.-2\,0.-3\,0$ &
$K12n320$ & $8^*-3.2\,1$ \\ \hline
\end{tabular}

\noindent
\begin{tabular}{|c|c|c|c|c|c|c|c|} \hline
$K12n323$ & $3.-2\,-1\,0.-2.2\,0$ & $K12n324$ &
$-3.-2\,-1\,0.-2.2\,0$ & $K12n325$ & $8^*2\,1.-3\,0$ \\ \hline
$K12n326$ & $8^*-2\,-1.2\,1$ & $K12n327$ & $2\,0.2\,1.-2\,0.2\,1\,0$
& $K12n330$ & $(-3,2\,1)\,(3,2+1)$ \\ \hline $K12n333$ &
$2\,4,-3,2\,1$ & $K12n334$ & $-2\,-4,3,2\,1$ & $K12n335$ &
$(-2\,-1\,-1,2\,1)\,(3,2)$ \\ \hline $K12n337$ &
$-2\,-3,2\,1\,1,2\,1$ & $K12n338$ & $2\,2\,1,4,-2\,-1$ & $K12n339$ &
$2\,2\,1,-4,2\,1$ \\ \hline $K12n341$ & $2\,2\,1,3\,1,-2\,-1$ &
$K12n342$ & $2\,2\,1,-3\,-1,2\,1$ & $K12n343$ &
$-2\,-2\,-1,3\,1,2\,1$ \\ \hline $K12n344$ & $(3\,1,-2\,-1)\,(3,2)$
& $K12n345$ & $(-3\,-1,2\,1)\,(3,2)$ & $K12n347$ &
$2\,3,3\,1,-2\,-1$ \\ \hline $K12n348$ & $-2\,-1:3\,0:2\,1\,0$ &
$K12n350$ & $2\,1:-3\,0:2\,1\,0$ & $K12n353$ &
$-2\,-1\,-1\,0:-2\,-1\,0:-2\,0$ \\ \hline $K12n356$ &
$8^*-2\,-1.2\,0.-2$ & $K12n357$ & $8^*-2\,-1.-2\,0.-2$ & $K12n358$ &
$2.-3\,0.-2\,-1.2\,0$ \\ \hline $K12n359$ & $3:2\,1:-2\,-1\,0$ &
$K12n360$ & $-2\,-1.2\,0.-2.2.-2\,0$ & $K12n361$ &
$8^*2\,0.2\,1:.-2$
\\ \hline $K12n363$ & $2\,1:2:-2\,-1\,-1\,0$ & $K12n364$ &
$.2.(-2\,-1,2\,1)$ & $K12n365$ & $.2.(2\,1,-2\,-1)$ \\ \hline
$K12n367$ & $2.-3.2\,1.2\,0$ & $K12n369$ &
$2\,0.2\,1.2\,0.-2\,-1\,0$ & $K12n372$ & $2\,1\,0.-2\,-1.2\,0.2\,0$
\\ \hline $K12n373$ & $2\,1\,0.2\,1.-2\,0.2\,0$ & $K12n375$ &
$2.2.-2.2\,0.-2\,-1$ & $K12n376$ & $2.2.-2.-2\,0.-2\,-1$ \\ \hline
$K12n378$ & $2.-2.2.2\,0.2\,1$ & $K12n379$ &
$2\,2\,1\,1,-2\,-1,2\,1$ & $K12n380$ & $-2\,-2\,-1\,-1,2\,1,2\,1$ \\
\hline $K12n381$ & $-2\,-1.-2.-2.2\,1\,0$ & $K12n383$ &
$-2\,-1.2.2.2\,1\,0$ & $K12n384$ & $2\,1.-2.2.2\,1\,0$ \\ \hline
$K12n385$ & $2\,1.2.-2.2\,1\,0$ & $K12n388$ &
$-2\,-2\,-1\,-1,3,2\,1$ & $K12n389$ & $2\,2\,1\,1,-3,2\,1$ \\ \hline
$K12n390$ & $2\,1\,1.2.-2.2\,0$ & $K12n391$ & $2\,1\,1.-2.2.2\,0$ &
$K12n392$ & $-3:2\,1:2\,1\,0$ \\ \hline $K12n393$ &
$8^*2.2\,0:-2\,-1\,0$ & $K12n395$ & $-3\,0:2\,1\,0:2\,1\,0$ &
$K12n396$ & $8^*2.-2\,0:2\,1\,0$ \\ \hline $K12n399$ &
$-2\,-1\,-1:-2\,-1\,0:-2\,0$ & $K12n400$ & $-2\,-1\,-1:2\,1\,0:2\,0$
& $K12n401$ & $2\,1\,1:2\,1\,0:-2\,0$ \\ \hline $K12n405$ &
$8^*2\,1:.-2\,-1\,0$ & $K12n406$ & $8^*2:.2\,0:.-2\,-1\,0$ &
$K12n407$ & $8^*2:.-2\,0:.-2\,-1\,0$ \\ \hline $K12n408$ &
$8^*2:.-2\,0:.2\,1\,0$ & $K12n409$ & $2\,0.-2\,-1.2\,0.3\,0$ &
$K12n410$ & $2\,0.2\,1.2\,0.-3\,0$ \\ \hline $K12n412$ &
$2\,0.3.2\,0.-2\,-1\,0$ & $K12n413$ & $2\,0.-3.2\,0.2\,1\,0$ &
$K12n415$ & $2\,1\,0.-3.2.2\,0$ \\ \hline $K12n416$ &
$-2\,-1:2\,1\,0:2\,1\,0$ & $K12n420$ & $8^*2\,1.-2\,0:.2$ &
$K12n421$ & $.2.(-2\,-1,2\,1)\,0$ \\ \hline $K12n422$ &
$.2.(2\,1,-2\,-1)\,0$ & $K12n424$ & $8^*2\,0.-2\,-1\,0:.2\,0$ &
$K12n427$ & $8^*-2\,0.2\,1\,0:.2\,0$ \\ \hline $K12n428$ &
$-2\,-1.2.2.2\,0.2$ & $K12n429$ & $8^*2:-2\,-1\,0:2\,0$ & $K12n431$
& $8^*2:2\,1\,0:-2\,0$
\\ \hline $K12n434$ & $8^*-2\,-2.2\,0$ & $K12n435$ &
$2.-2\,-1.2\,0.2\,1\,0$ & $K12n440$ & $8^*2\,1.-2:.2$ \\ \hline
$K12n441$ & $8^*-3::2\,1$ & $K12n443$ & $8^*-2\,-1::-3\,0$ &
$K12n445$ & $8^*2\,1::-3\,0$ \\ \hline $K12n447$ &
$.2\,1.-2\,-1.2.2\,0$ &
$K12n448$ & $8^*2\,0.2\,1:-2\,0$ & $K12n450$ & $8^*2.2\,1.-2\,0$ \\
\hline $K12n451$ & $8^*2:.-2\,-1\,0:.-2\,0$ & $K12n452$ &
$3.-2\,0.-2.2\,1\,0$ & $K12n453$ & $8^*2:.2\,1\,0:.-2\,0$ \\ \hline
$K12n454$ & $(4,-2\,-1)\,(3,2)$ & $K12n455$ & $(-4,2\,1)\,(3,2)$ &
$K12n458$ & $8^*-2\,-1\,-2\,0$ \\ \hline $K12n459$ &
$8^*2\,1\,1\,0::-2\,0$ & $K12n460$ & $8^*2\,2.-2\,0$ & $K12n461$ &
$102^*:-2\,-1\,0$ \\ \hline $K12n462$ & $8^*-2\,0.2:-2\,-1\,0$ &
$K12n463$ & $8^*-2\,0.2:2\,1\,0$ & $K12n465$ & $102^*2\,0:.-2$ \\
\hline $K12n466$ & $4\,1,3\,1,-2\,-1$ & $K12n467$ &
$4\,1,-3\,-1,2\,1$ & $K12n468$ & $-4\,-1,3\,1,2\,1$ \\ \hline
$K12n469$ & $3\,1\,1,3\,1,-2\,-1$ & $K12n470$ &
$3\,1\,1,-3\,-1,2\,1$ & $K12n471$ & $-3\,-1\,-1,3\,1,2\,1$ \\ \hline
$K12n474$ & $5,3\,1,-2\,-1$ & $K12n476$ & $-5,3\,1,2\,1$ & $K12n477$
& $3\,2,3\,1,-2\,-1$ \\ \hline $K12n478$ & $3\,2,-3\,-1,2\,1$ &
$K12n479$ & $-3\,-2,3\,1,2\,1$ & $K12n480$ & $-3:2\,1:2\,1$ \\
\hline $K12n481$ & $-2.-2\,0.-2.2\,2\,0$ & $K12n482$ &
$-2\,-1\,-1:2\,1:2\,0$ & $K12n484$ & $8^*-2\,-1.2:.2$ \\ \hline
$K12n485$ & $-2\,-1.-2.-2.2.2$ & $K12n486$ &
$2\,1\,0.-2.-2\,0.-3\,0$ & $K12n489$ & $2\,1\,0.2.2\,0.-3\,0$ \\
\hline $K12n490$ & $-3.-2\,0.-2.2\,1\,0$ & $K12n491$ &
$8^*2:.-2\,-1\,0:.2\,0$ & $K12n492$ & $2\,1\,0.-3.2\,0.2\,0$ \\
\hline $K12n493$ & $2\,1\,0.2.-3.2\,0$ & $K12n496$ &
$(-3,2\,1)\,(2\,1,2\,1)$ & $K12n497$ & $2\,1\,2\,1,-3,2\,1$ \\
\hline $K12n498$ & $8^*2\,0.2:-2\,-1\,0$ & $K12n500$ &
$.-2\,-1.4.2\,0$ & $K12n501$ & $.-2\,-1.3\,1.2\,0$ \\ \hline
$K12n504$ & $8^*-2\,-2\,0:.2\,0$ & $K12n505$ & $8^*2\,2\,0:.-2\,0$ &
$K12n506$ & $8^*2\,0.2:.-2\,-1\,0$ \\ \hline $K12n507$ &
$9^*2:-2\,-1\,0$ & $K12n508$ & $(-2\,-1,2\,1)\,(2\,1,2\,1)$ &
$K12n509$ & $-3:-2\,-1\,0:-2\,-1\,0$ \\ \hline $K12n510$ &
$-3:2\,1\,0:2\,1\,0$ & $K12n511$ & $2\,1\,0.2.2\,0.-2\,-1\,0$ &
$K12n512$ & $8^*2\,0.2\,0:.-2\,-1\,0$ \\ \hline $K12n513$ &
$8^*-2\,0.2\,0:.2\,1\,0$ & $K12n514$ & $-2\,-1.2.2.2.2$ & $K12n515$
& $2\,1.-2.2.2.2$ \\ \hline $K12n516$ & $2\,1.2.2.-2.2$ & $K12n517$
& $2.-2.2\,1.2\,0.2$ & $K12n520$ & $2.2\,1.-3.2\,0$ \\ \hline
$K12n521$ & $8^*2.-2\,0.2\,1\,0$ & $K12n522$ & $4\,1,4,-2\,-1$ &
$K12n524$ & $-3\,0:-2\,-1\,0:-2\,-1\,0$ \\ \hline $K12n525$ &
$2\,1:2\,1:-3\,0$ & $K12n527$ & $2.2\,1.-2.2.2\,0$ & $K12n529$ &
$2.-2.2\,1.2\,1\,0$ \\ \hline $K12n530$ & $8^*2\,0.-2\,-1\,0:2\,0$ &
$K12n531$ & $2\,1.-2.2.2\,0.2$ & $K12n532$ & $9^*2.-2\,-1$ \\ \hline
$K12n533$ & $8^*-2\,-1.2:.2\,0$ & $K12n534$ & $8^*2.2\,0.-2\,-1\,0$
& $K12n536$ & $8^*-2.-2\,0.-2\,-1\,0$ \\ \hline $K12n537$ &
$102^*-2\,-1\,0$ & $K12n538$ & $9^*-2\,-1:.2\,0$ & $K12n539$ &
$9^*2\,1:.-2\,0$ \\ \hline $K12n540$ & $8^*2\,0.-2\,-1\,0:.2$ &
$K12n541$ & $2.-2\,-1.3.2\,0$ & $K12n542$ & $-2.-2\,-1.-2.2.2\,0$ \\
\hline $K12n543$ & $2.2.-2\,-1.2.2$ & $K12n544$ & $2.-2.-2\,-1.-2.2$
& $K12n545$ & $2.-2.2\,1.2.2$ \\ \hline $K12n546$ &
$.2.-2\,-1.2.3\,0$ & $K12n547$ & $.3.2.-2\,-1.2$ & $K12n548$ &
$8^*2.-2\,-1.2$ \\ \hline $K12n551$ &
$-2\,-1\,-1\,-1\,0:-2\,0:-2\,0$ & $K12n553$ & $3,2\,1,2\,1,-2\,-1$ &
$K12n555$ & $-3,2\,1,2\,1,2\,1$
\\ \hline $K12n556$ & $3,2\,1,-2\,-1,2\,1$ & $K12n557$ &
$-2\,-2\,0:-2\,-1\,0:-2\,0$ & $K12n559$ & $8^*3\,0:.-2:.2\,0$ \\
\hline $K12n560$ & $8^*2\,1:-3\,0$ & $K12n561$ & $8^*3:.-2:.2\,0$ &
$K12n562$ & $.2\,1.-2.4\,0$ \\ \hline $K12n563$ & $3:2\,2\,0:-2\,0$
& $K12n564$ & $2\,1\,2\,1,3,-2\,-1$ & $K12n565$ & $-4\,-2,3,2\,1$ \\
\hline $K12n566$ & $-3\,-1\,-2,3,2\,1$ & $K12n567$ &
$9^*-2\,0.2\,1\,0$ & $K12n568$ & $9^*2\,0.2::-2\,0$ \\ \hline
$K12n569$ & $-2\,-2\,0:-3\,0:-2\,0$ & $K12n570$ & $5\,1,2\,1,-2\,-1$
& $K12n571$ & $-5\,-1,2\,1,2\,1$ \\ \hline $K12n572$ &
$3\,2\,1,2\,1,-2\,-1$ & $K12n573$ & $-3\,-2\,-1,2\,1,2\,1$ &
$K12n576$ & $3\,2\,1,3,-2\,-1$ \\ \hline $K12n578$ &
$3\,2\,1,-3,2\,1$ & $K12n580$ & $-3:-2\,-2\,0:-2\,0$ & $K12n581$ &
$5\,1,3,-2\,-1$ \\ \hline $K12n583$ & $5\,1,-3,2\,1$ & $K12n584$ &
$2\,1.-2.3.2$ & $K12n585$ & $-2\,-1.3.2.2$ \\ \hline $K12n586$ &
$2\,1.3.-2.2$ & $K12n587$ & $2.-2.2\,1\,1.2\,0$ & $K12n588$ &
$9^*-2\,-1\,0.2\,0$ \\ \hline $K12n589$ & $9^*2\,1\,0:-2\,0$ &
$K12n590$ & $3\,0.-2\,-1.2\,0.2\,0$ & $K12n592$ &
$3\,0.2\,1.-2\,0.2\,0$ \\ \hline $K12n593$ & $3.-2\,-1.2\,0.2\,0$ &
$K12n595$ & $3.2\,1.-2\,0.2\,0$ & $K12n596$ & $-2.2.2\,0.2\,1\,1$ \\
\hline $K12n597$ & $9^*-2\,-1\,0:2\,0$ & $K12n598$ & $8^*-2\,-3$ &
$K12n599$ & $8^*-2\,-1.2:2$ \\ \hline $K12n600$ &
$(3,3)\,(3,-2\,-1)$ & $K12n602$ & $(3,3)\,(-3,2\,1)$ & $K12n604$ &
$(3,3)\,(-2\,-1,2\,1)$ \\ \hline $K12n606$ & $8^*-2\,-1\,0:2:2\,0$ &
$K12n607$ & $-2\,-1.3.2\,0.2$ & $K12n608$ & $-2\,-1.2.3.2$ \\ \hline
$K12n609$ & $8^*2.-2\,-1:2$ & $K12n610$ & $3\,0.2.2\,0.-2\,-1\,0$ &
$K12n611$ & $3.2.2\,0.-2\,-1\,0$ \\ \hline $K12n612$ &
$-2\,-2\,-1\,0:-2\,0:-2\,0$ & $K12n613$ & $8^*2\,1\,1:-2\,0$ &
$K12n614$ & $8^*-2\,-1\,-1\,0:2\,0$ \\ \hline $K12n615$ &
$8^*2\,1\,1\,0:-2\,0$ & $K12n616$ & $9^*-2\,-1\,0.2$ & $K12n618$ &
$8^*2\,1\,0.-2\,0:2\,0$ \\ \hline $K12n619$ &
$8^*2\,1\,0:2\,0:-2\,0$
& $K12n620$ & $9^*2\,1\,0.-2$ & $K12n621$ & $8^*2\,1:2\,0:-2\,0$ \\
\hline $K12n622$ & $8^*2\,1:-2\,0:2\,0$ & $K12n623$ &
$8^*-2\,-1\,0:.2:.2\,0$ & $K12n625$ & $8^*2\,1\,0:.2:.-2\,0$ \\
\hline $K12n626$ & $9^*.-2:.2\,1\,0$ & $K12n627$ &
$2\,1\,0.2.-2.2.2\,0$ & $K12n628$ & $9^*.-2\,-1:.2$ \\ \hline
$K12n631$ & $8^*2\,1\,0:-2\,0:2\,0$ & $K12n632$ & $9^*.2\,1:.-2$ &
$K12n633$ & $8^*-2\,-1\,0.2\,0:2\,0$ \\ \hline

\end{tabular}

\noindent
\begin{tabular}{|c|c|c|c|c|c|c|c|} \hline
$K12n635$ & $9^*2\,1.-2$ & $K12n636$ & $8^*2\,2\,0:-2\,0$ &
$K12n637$ & $9^*.-2\,-1\,-1$ \\ \hline $K12n641$ &
$3\,1:2\,1\,0:-2\,0$ &
$K12n643$ & $4:-2\,-1\,0:2\,0$ & $K12n645$ & $4:2\,1\,0:-2\,0$ \\
\hline $K12n646$ & $8^*2\,0.-2\,-1\,0:2$ & $K12n649$ &
$8^*-2\,-1\,0.2\,0.2\,0$ & $K12n651$ & $8^*2\,1\,0.2\,0.-2\,0$ \\
\hline $K12n652$ & $8^*2\,1\,0.-2\,0.2\,0$ & $K12n653$ &
$-2\,-1.2.2.3\,0$ & $K12n654$ & $3.2.2.-2\,-1$ \\ \hline $K12n656$ &
$3.-2.2.2\,1$ & $K12n657$ & $2.-2\,-1.2.3\,0$ & $K12n658$ &
$8^*-2\,-1\,0.2\,0:2$ \\ \hline $K12n659$ & $3.-2\,-1.2.2$ &
$K12n661$ & $2:4\,0:-2\,-1\,0$ & $K12n662$ & $3.2.-2\,-1.2$ \\
\hline $K12n663$ & $9^*2\,0::-2\,-1$ & $K12n664$ & $9^*.-2\,-2$ &
$K12n666$ & $8^*-5$ \\ \hline $K12n667$ & $8^*-3\,-2$ & $K12n668$ &
$8^*-3\,-2\,0$ & $K12n669$ & $8^*-5\,0$ \\ \hline $K12n671$ &
$(3,2),-2,(3,2)$ & $K12n674$ & $(4\,1,-2)\,(3,2)$ & $K12n675$ &
$(-4\,-1,2)\,(3,2)$ \\ \hline $K12n677$ & $(3\,1\,1,-2)\,(3,2)$ &
$K12n678$ & $(-3\,-1\,-1,2)\,(3,2)$ & $K12n682$ & $(3,2),-2,(2,3)$
\\ \hline $K12n683$ & $2.-2.3.2.2$ & $K12n684$ & $8^*-3:.3$ &
$K12n685$ & $8^*-3:.3\,0$ \\ \hline $K12n686$ & $8^*3:.-3\,0$ &
$K12n687$ & $8^*-3\,0:.3\,0$ & $K12n695$ & $2:(2,3)\,0:-2\,0$ \\
\hline $K12n698$ & $8^*2:.-4\,0$ & $K12n699$ & $8^*2:.-3\,-1\,0$ &
$K12n700$ & $8^*3\,1:.-2\,0$ \\ \hline $K12n701$ & $8^*4:.-2\,0$ &
$K12n702$ & $102^*.-2.2$ & $K12n703$ & $102^*.-2:2$ \\ \hline
$K12n704$ & $9^*3\,0:.-2\,0$ & $K12n705$ & $101^*2:::.-2\,0$ &
$K12n706$ & $101^*-2\,0::.-2\,0$ \\ \hline $K12n707$ & $2.2.2.-2.3$
& $K12n708$ & $-2.2.3.2.2$ & $K12n709$ & $8^*-3:3$ \\ \hline
$K12n710$ & $8^*-3:3\,0$ & $K12n711$ & $8^*3:-3\,0$ & $K12n712$ &
$9^*3:.-2\,0$
\\ \hline $K12n713$ & $8^*3\,0:2:-2\,0$ & $K12n714$ & $9^*-3\,0.2\,0$
& $K12n715$ & $8^*3\,0.2:-2\,0$ \\ \hline $K12n716$ & $9^*-3.2\,0$ &
$K12n717$ & $8^*-2.-2\,0.-2\,0:2$ & $K12n718$ & $8^*2.3.-2\,0$ \\
\hline $K12n719$ & $2.-2.-2.-2.3\,0$ & $K12n720$ &
$8^*2\,0.-2\,0.2:2\,0$ & $K12n721$ & $-4\,-1,4\,1,2$ \\ \hline
$K12n722$ & $4\,1,4\,1,-2$ & $K12n723$ & $-4\,-1,3\,1\,1,2$ &
$K12n724$ & $4\,1,3\,1\,1,-2$ \\ \hline $K12n726$ &
$3\,1\,1,3\,1\,1,-2$ & $K12n727$ & $8^*-3\,0:3\,0$ & $K12n728$ &
$102^*2\,0::-2$ \\ \hline $K12n729$ & $102^*-2\,0::-2$ & $K12n730$ &
$8^*-2.2.-2\,0:2\,0$ & $K12n731$ & $102^*-2\,0.2$ \\ \hline
$K12n732$ & $101^*2::-2\,0$ & $K12n733$ & $8^*-3\,0.3\,0$ &
$K12n734$ & $8^*2\,0.3\,0:.-2\,0$ \\ \hline $K12n735$ &
$8^*2\,0.3:.-2\,0$ &
$K12n736$ & $102^*-2\,0::2$ & $K12n737$ & $8^*-3\,-1\,0::2\,0$ \\
\hline $K12n738$ & $8^*2.2.2\,0.-2\,0$ & $K12n740$ &
$8^*2.2.-2\,0.2\,0$ & $K12n741$ & $8^*2:.3\,0:.-2\,0$ \\ \hline
$K12n742$ & $9^*-3::2\,0$ & $K12n743$ & $9^*:3\,0.-2\,0$ & $K12n744$
& $2.3.-2.2\,0.2\,0$ \\ \hline $K12n745$ & $8^*2.3\,0.-2\,0$ &
$K12n746$ & $9^*-3\,0::2\,0$ & $K12n747$ & $2.-2.-2.-2.3$ \\ \hline
$K12n748$ & $2.2.2.2.-3$ & $K12n751$ & $2.2.-2.2.3$ & $K12n752$ &
$2.2.-2.2.3\,0$ \\ \hline $K12n753$ & $8^*2\,0.2\,0.2:-2\,0$ &
$K12n754$ & $8^*2\,0.-2\,0:.3\,0$ & $K12n755$ &
$8^*2\,0.2\,0:.-3\,0$
\\ \hline $K12n756$ & $8^*2\,0.-2\,0:.3$ & $K12n757$ &
$8^*2\,0.-2\,0:2\,0:2$ & $K12n758$ & $8^*2\,0.2\,0:-2\,0:2$ \\
\hline $K12n759$ & $8^*3\,0.-2\,0:.2$ & $K12n760$ & $9^*:-3\,0.2\,0$
& $K12n761$ & $8^*2.2\,0:-2\,0:2\,0$ \\ \hline $K12n762$ &
$8^*2.2\,0:2\,0:-2\,0$ & $K12n763$ & $9^*2.-2:.2\,0$ & $K12n765$ &
$9^*2.-2.2\,0$ \\ \hline $K12n766$ & $9^*2.2.-2\,0$ & $K12n767$ &
$2:-4\,0:3\,0$ & $K12n769$ & $8^*.2:2\,0.2\,0:-2\,0$ \\ \hline
$K12n770$ & $9^*2\,0:.2.-2$ & $K12n771$ & $8^*.2:2\,0.-2\,0:2\,0$ &
$K12n772$ & $8^*.2:-2\,0.2\,0:2\,0$ \\ \hline $K12n773$ &
$3:2:-4\,0$
& $K12n774$ & $3:2:-3\,-1\,0$ & $K12n775$ & $102^*-2\,0:2\,0$ \\
\hline $K12n776$ & $102^*2\,0:-2\,0$ & $K12n777$ &
$8^*3\,0:.2:.-2\,0$ & $K12n778$ & $2.-2.3.2\,0.2$ \\ \hline
$K12n779$ & $2.2.-3.2.2\,0$
& $K12n780$ & $8^*2.2.2\,0:-2\,0$ & $K12n781$ & $8^*-2.2.2\,0:2\,0$ \\
\hline $K12n782$ & $8^*-2.2.2:2\,0$ & $K12n783$ &
$8^*-3\,0.2\,0.2\,0$ & $K12n784$ & $9^*2\,0:.2\,0.-2$ \\ \hline
$K12n785$ & $8^*3:.2:.-2\,0$ & $K12n786$ & $8^*-3\,0.2\,0:.2$ &
$K12n787$ & $3.2.-2.-2\,0.-2$ \\ \hline $K12n788$ &
$-2.-2\,0.-2.3.2\,0$ & $K12n789$ & $8^*2\,0.2\,0.-2\,0:2\,0$ &
$K12n790$ & $8^*2\,0.-2\,0.2\,0:2\,0$ \\ \hline $K12n791$ &
$8^*-2\,0.2\,0.2\,0:2\,0$ & $K12n792$ & $2.-4.-2\,0.-2$ & $K12n793$
& $-2.4.2\,0.2$ \\ \hline $K12n794$ & $2.-3\,-1.-2\,0.-2$ &
$K12n795$ & $-2.3\,1.2\,0.2$ & $K12n796$ & $8^*2.3:-2\,0$ \\ \hline
$K12n797$ & $8^*2.3\,0:-2\,0$ & $K12n798$ & $-2.-2.-2.2.2\,0.2\,0$ &
$K12n799$ & $9^*2.-2:2\,0$ \\ \hline $K12n800$ & $9^*2.2:-2\,0$ &
$K12n802$ & $8^*2\,0.2:-2\,0.2\,0$ & $K12n803$ &
$8^*2\,0.2:2\,0.-2\,0$ \\ \hline $K12n804$ & $3\,1.-2.2.2\,0$ &
$K12n805$ & $4.-2.2.2\,0$ & $K12n806$ & $8^*-2\,0:-2\,0:-2\,0:-2\,0$
\\ \hline $K12n807$ & $8^*-2\,0:-2\,0:-2\,0:2\,0$ & $K12n811$ &
$9^*.-2:-2\,0.-2$ & $K12n813$ & $9^*.2:-2\,0.2$ \\ \hline $K12n814$
& $8^*2\,0:-3\,0:2\,0$
& $K12n816$ & $9^*.-2:2\,0.2$ & $K12n817$ & $-3\,-1.2.2\,0.2$ \\
\hline $K12n818$ & $-4.2.2\,0.2$ & $K12n819$ & $8^*2\,0.2:-2\,0:2$ &
$K12n820$ & $3\,0.2.2\,0.-3\,0$ \\ \hline $K12n821$ & $2:5\,0:-2\,0$
& $K12n822$ & $-3\,0.2.2\,0.3\,0$ & $K12n823$ & $-3\,-2:2\,0:2\,0$
\\ \hline $K12n824$ & $2:3\,2\,0:-2\,0$ & $K12n825$ & $-2.2.3.3\,0$
& $K12n826$ & $8^*.2:2\,0.2:-2\,0$ \\ \hline $K12n827$ &
$9^*2:2\,0.-2\,0$ & $K12n828$ & $9^*-2\,0.2\,0::2\,0$ & $K12n829$ &
$2.-2.-3.-2.2$ \\ \hline $K12n831$ & $2.2.-3.2.2$ & $K12n832$ &
$2\,0.2.-3.-2\,0.-2$ & $K12n833$ & $8^*2\,0.2.2:-2\,0$ \\ \hline
$K12n834$ & $8^*2\,0.-2\,0:2\,0.2\,0$ & $K12n836$ &
$8^*-2\,0.2\,0:2\,0.2\,0$ & $K12n837$ & $1212^*-1.-1.-1.-1.-1.-1$ \\
\hline $K12n839$ & $111^*:.-2\,0$ & $K12n840$ & $101^*-2\,0::.2\,0$
& $K12n841$ & $8^*2.-2\,0.2\,0.2\,0$ \\ \hline $K12n842$ &
$111^*:::-2\,0$ & $K12n843$ & $9^*-2\,0:.-2\,0:.-2\,0$ & $K12n844$ &
$8^*2.2\,0.2:.-2\,0$ \\ \hline $K12n845$ & $9^*-2\,0::3$ & $K12n846$
& $8^*4:-2\,0$ & $K12n847$ & $9^*.-4\,0$ \\ \hline $K12n848$ &
$8^*3\,0.2\,0.-2\,0$ & $K12n849$ & $9^*2\,0::-3$ & $K12n852$ &
$8^*3\,0.-2\,0.2\,0$ \\ \hline $K12n853$ & $8^*2.2:-2\,0.2\,0$ &
$K12n854$ & $9^*.-2:2.2\,0$ & $K12n855$ & $3\,2:2:-2\,0$ \\ \hline
$K12n856$ & $5:2:-2\,0$ & $K12n857$ & $102^*-2\,0:.2$ & $K12n858$ &
$8^*-2\,0.2\,0.2\,0:2$ \\ \hline $K12n859$ & $9^*2.-2::2\,0$ &
$K12n860$ & $8^*2\,0.-3\,0:2$ & $K12n861$ & $9^*-2\,0.2::2\,0$ \\
\hline $K12n862$ & $9^*2\,0.-2::2\,0$ & $K12n863$ &
$9^*2\,0.2\,0:.-2$ & $K12n864$ & $8^*2\,0.-2\,0.2\,0:2$ \\ \hline
$K12n865$ & $3\,0.2.-2.2.2\,0$ & $K12n866$ & $8^*2.2\,0.-2\,0.2\,0$
& $K12n867$ & $8^*-3\,0.2\,0:2$ \\ \hline $K12n869$ &
$9^*.-2:.2:.2\,0$ & $K12n871$ & $9^*2\,0.-2:.2$ & $K12n872$ &
$9^*2\,0.2:.-2$ \\ \hline $K12n873$ & $8^*2\,0.-2\,0.-2\,0.2\,0$ &
$K12n874$ & $8^*2\,0.-2\,0.2\,0.2\,0$ & $K12n875$ & $2.-2.3\,1.2\,0$
\\ \hline $K12n876$ & $2.-2.4.2\,0$ &
$K12n877$ & $102^*-2\,0::.2\,0$ & $K12n878$ & $9^*2.-2:::.2\,0$ \\
\hline $K12n879$ & $102^*:-2\,0.2\,0$ & $K12n880$ & $102^*.2:-2\,0$
& $K12n881$ & $2.-2.-2\,0.-2.2.2\,0$ \\ \hline $K12n882$ &
$2\,0.3.2\,0.-3\,0$ & $K12n883$ & $9^*-2\,0.3\,0$ & $K12n884$ &
$9^*-2\,0.3$ \\ \hline $K12n885$ & $9^*2\,0.-3$ & $K12n886$ &
$9^*.-3:.2\,0$ & $K12n887$ & $8^*2.-3.2$ \\ \hline

\end{tabular}

\normalsize


\bigskip
\bigskip

\footnotesize

\noindent THE MATHEMATICAL INSTITUTE, KNEZ MIHAILOVA  36, P.O.BOX
367, \\ 11001 BELGRADE, SERBIA

\medskip

\noindent {\it E-mail address:} $\mathrm{sjablan@gmail.com}$


\begin{thebibliography}{99}

\bibitem [BN]{BN} Bar-Natan, D. (2002) On Khovanov's categorification
of the Jones polynomial, {\it Alg. Geom. Top.} {\bf 2}, 337--370.
\medskip
\bibitem [BiWi]{BiWi} Birman, J.~S. and Williams, R.~F. (1983) Knotted
Periodical Orbits in Dynamical Systems-I: Lorenz's Equations,
Topology {\bf 22}, 1, 47--82.
\medskip
\bibitem [Cau]{Cau} Caudron, A. (1982) Classification des n\oe uds et des enlancements,
Public. Math. d'Orsay 82. Univ. Paris Sud, Dept. Math., Orsay.
\medskip
\bibitem [ChKo1]{ChKo1} Champanerkar, A. and Kofman, I. (2008) Twisting
quasi-alternating links, \\ arXiv:math/0712.2990v2 [math.GT]
\medskip
\bibitem [ChKo2]{ChKo2} Champanerkar, A. and Kofman, I. (2005) On the Mahler measure of Jones polynomials under twisting
   \emph{ Alg. Geom. Top.}  \textbf{5}, 1–-22
\medskip
\bibitem [Con]{Con} Conway, J. (1970) An enumeration of knots and links and some of
their related properties, in {\it Computational Problems in Abstract
Algebra}, Proc. Conf. Oxford 1967 (Ed. J. Leech), 329--358, Pergamon
Press, New York.
\medskip
\bibitem [Don]{Don}  Donaldson S.~K. (1987) The orientation of Yang-Mills moduli spaces and 4-manifold topology. J. Diferential
Geom., {\bf 26}, 3, 397–-428.
\medskip
\bibitem [Gh]{Gh}Ghys, E. (2006) Knots and Dynamics, preprint, to
appear in Proc. ICM-2006, Madrid.
\medskip
\bibitem [Gr]{Gr}Greene, J. (2009) Homologically thin,
non-quasialternating links, arXiv:math/0906.2222v1 [math.GT]
\medskip
\bibitem [GhLe]{GhLe}Ghys, E. and Leys, J. (2006) Lorenz and Modular
Flows: A Visual Introduction, AMS Feature Column, Nov. 2006,
http://www.ams.org/featurecolumn/archive/lorenz.html
\medskip
\bibitem [HosThi]{HosThi} Hoste, J. and Thistlethwaite, M. {\it
Knotscape}, http://www.math.utk.edu/~morwen/knotscape.html
\medskip
\bibitem [Ja]{Ja}Jablan, S.~V. (2008) Adequacy of link families,
arXiv:math/08011.0081v1 [math.GT]
\medskip
\bibitem [JaSa]{JaSa} Jablan, S.~V., Sazdanovi\' c, R. (2007) {\it
LinKnot- Knot Theory by Computer}. World Scientific, New Jersey,
London, Singapore.
\medskip
\bibitem [JaSa1]{JaSa1} Jablan, S.~V., Sazdanovi\' c, R. (2008)
Quasi-alternating links and odd homology: computations and
conjectures, arXiv:math/0901.0075v2 [math.GT]
\medskip
\bibitem [Kh1]{Kh1}Khovanov, M. (1999) A categorification of the Jones polynomial,
arXiv:math/9908171v2 [math.QA]
\medskip
\bibitem [Kh2]{Kh2} Khovanov, M. (2002) Patterns in knot homology I,
Experimental Mathematics 12, 3 (2003), 365-374 arXiv:math/0201306v1
[math.QA]
\medskip
\bibitem [Le]{Le} Lee, E.S. (2002) The support of the Khovanov's
invariant for alternating knots, arXiv:math/0201105v1 [math.GT]
\medskip
\bibitem [Lo]{Lo}Lorenz, E.N. (1963) Deterministic non-periodic
flow, J. Atmospheric Science {\bf 20}, 130--141.
\medskip
\bibitem [MaOz]{MaOz} Manolescu, C. and Ozsv\' ath, P. (2008) On the
Khovanov and knot Floer homologies of quasi-alternating links,
arXiv:math/0708.3249v2 [math.GT]
\medskip
\bibitem [OzRaSz]{OzRaSz}Ozsv\' ath, P., Rasmussen, J. and Szab\' o, Z. (2007)
Odd Khovanov homology, arXiv:math/0710.4300v1 [math.QA]).
\medskip
\bibitem [OzSz]{OzSz}Ozsv\' ath, P. and Szab\' o, Z. (2005) On the
Heegaard Floer homology of branched double-covers, {\it Adv. Math.}
194 {\bf 1}, 1--33 (arXiv:math/0309170v1 [math.GT]).
\medskip
\bibitem [Sh1]{Sh1}Shumakovitch, A. (2008) KhoHo,
http://www.geometrie.ch/KhoHo/
\medskip
\bibitem [Sh2]{Sh2} Shumakovitch, A. (2008) Private communication.
\medskip
\bibitem [St]{St}Sto\v si\' c, M. (2006) Homological thickness and
stability of torus knots, arXiv:math/0511532v2 [math.GT]
\medskip
\bibitem [Wa]{Wa} Watson, M. (2008) Surgery obstructions from Khovanov homology, arXiv:0807.1341v2 [math.GT]
\medskip
\bibitem [Wi]{Wi}Widmer, T. (2008) Quasi-alternating Montesinos links, arXiv:math/0811.0270v1 [math.GT]

\end{thebibliography}
\end{document}